\documentstyle[12pt]{amsart}
\def\ds{\displaystyle}

\newsymbol\Subset 1362

\newcommand{\R}{{\Bbb R}}

\def\tl{\widetilde}
\def\ol{\overline}

\def\C{{\cal C}}
\def\E{{\cal  E}}

\def\K{{\cal K}}

\def\H{{\cal H}}

\newtheorem{thm}{Theorem}[subsection]  
\newtheorem{lem}[thm]{Lemma}	       
\newtheorem{crlr}[thm]{Corollary}      
\newtheorem{defin}[thm]{Definition}    
\newtheorem{rem}[thm]{Remark}	       

\makeatletter
\def\theequation{\thesection.\@arabic\c@equation}
\def\thethm{\thesection.\@arabic\c@thm}
\def\thelem{\thesection.\@arabic\c@thm}
\def\thecrlr{\thesection.\@arabic\c@thm}
\def\theprp{\thesection.\@arabic\c@thm}
\def\therem{\thesection.\@arabic\c@thm}
\makeatother

\begin{document}

\baselineskip=18pt

\title[Singular Integrals and Commutators  in Generalized Morrey Spaces]{Singular
Integrals and Commutators  in Generalized Morrey Spaces}

\author [L.G. Softova]{Lubomira G. Softova}
\address{ Bulgarian Academy of Sciences,
Institute of Mathematics and Informatics,
Acad. G. Bonchev Str. bl. 8, 1113 Sofia, Bulgaria}
\curraddr{Department of Mathematics, University of Bari, 4~E. Orabona~Str.,
 70~125~Bari,	Italy}
\email{luba{@@}pascal.dm.uniba.it}

\begin{abstract}
The purpose of this paper is to study singular integrals whose kernels
$k(x;\xi)$  are variable, i.e. they  depend on some parameter $x\in \R^n$	 and in
$\xi\in\R^n\setminus\{0\}$ satisfy mixed homogeneity condition of the form
$k(x;\mu^{\alpha_1}\xi_1,\ldots,\mu^{\alpha_n}\xi_n)=
\mu^{-\sum_{i=1}^n \alpha_i}k(x;\xi)$	with positive real numbers
$\alpha_i\geq 1$ and  $\mu>0.$	The continuity of these operators in
$L^p(\R^n)$ is
well studied by Fabes and Rivi\`ere. Our goal is to extend their results in generalized
Morrey spaces with a weight satisfying suitable dabbling and integral conditions. A special attention is paid also of the
commutators of the kernel with functions of bounded and vanishing mean
oscillation.
\end{abstract}

\maketitle

\section{Introduction}\label{s1}
\setcounter{equation}{0}
\setcounter{thm}{0}

We  consider  the following integral operators
\begin{equation}\label{Kf}
\K f(x):=  P.V. \int_{\R^n} k(x;x-y)f(y)dy
\end{equation}
and its commutators with essentially bounded functions
\begin{align}
\nonumber
\C[a,k]f(x):= &\  P.V.\int_{\R^n} k(x;x-y)[a(y)-a(x)]f(y)dy\\
\label{Cakf}
   = &\  \K(af)(x)-a(x)\K f(x).
\end{align}
The generating kernel $k(x;\xi) : \R^n\times\R^n\setminus\{0\}\to \R$
is variable, i.e. it depends on some parameter $x$   and possesses ``good''
properties with respect to the second variable $\xi.$  This class
 of kernels is firstly studied by Fabes and Rivi\`{e}re in
\cite{FR}.  They generalize the classical kernels of Calder\'on and Zygmund
$k(\xi)=\Omega(\xi)/|\xi|^{n}$	having homogeneity of degree $-n$	 and those studied by
Jones in \cite{BJ}   and satisfying homogeneity property of the form
$k(\lambda \xi,\lambda^m \tau)=\lambda^{-n-m}k(\xi,\tau),$  $\xi\in \R^n,$
	$\tau\in(0,\infty),$  $m\geq
1.$  Introducing a new metric $\rho,$	Fabes and Rivi\`ere study  \eqref{Kf}
 in $L^p(\R^n),$ where	$\R^n$ is endowed with the topology induced by $\rho$
	and defined by ellipsoids. Thus,
the unite sphere with respect to $\rho$ coincides with the unite sphere
$\Sigma_n$  with respect to the Euclidean metric. This fact allows to impose on the kernel
$k$ the Calder\'on-Zygmund conditions on the unite sphere, in spite of the lack of
"symmetry" of $k$  with respect to the variables $\xi_i,$  $i=1,\ldots,n.$  Let we
note, that the standard parabolic metric $\tl\rho=\sup\{|x|,\sqrt t\},$
$x\in\R^n,t\in(0,\infty),$ for instance,  does not  permit to define the
mentioned above conditions on kernels having homogeneity of parabolic type.
Using the Fourier transform in $L^2(\R^n)$ and the Marcinkiewicz interpolation theorem,
 Fabes and Rivi\`ere  obtained	that the  integral operators \eqref{Kf} are continuous in $L^p(\R^n),$ $p\in(1,\infty).$

In  the present work we study the continuity of these operators in the generalized Morrey spaces  $L^{p,\omega}(\R^n)$
   where the  function $\omega$   satisfies
suitable conditions.

A special attention is paid also to the commutators $\C[a,k]$	 of the
kernel $k$  and functions $a$  having bounded or vanishing mean oscillation.	In this case we
impose of the results of Coifman-Rochberg-Weiss (\cite{CRW})  and Bramanti-Cerutti (\cite{BC})	treating
continuity in $L^p(\R^n)$   of commutators with constant kernels.

The technique we  used is the one elaborated  by Calder\'on and Zygmund and
consisting  of expansion of the kernel into spherical harmonics  and restricting the considerations on integral operators with constant kernels.

\section{Definitions and preliminary results}\label{s2}

\setcounter{equation}{0}
\setcounter{thm}{0}

Let  $\alpha_1,\ldots,\alpha_n$ be real numbers, $\alpha_i\geq 1$ and define
$\alpha=\sum_{i=1}^n\alpha_i.$
Following Fabes and Rivi\`ere  (\cite{FR}),
the function
$
F(x,\rho)=\sum_{i=1}^{n} x_i^2 \rho^{-2\alpha_i},
$
considered for any fixed $x\in\R^n,$ is a decreasing one with respect to
$\rho>0$ and  the equation $F(x,\rho)=1$  is unique solvable in $\rho(x).$
 It is a
simple matter to check that $\rho(x-y)$ defines a distance between any two
points $x,y\in\R^n.$ Thus $\R^n,$ endowed with the metric $\rho$ results a
homogeneous metric space (\cite[Remark~1]{FR}, \cite{BC}).
The balls with respect to $\rho(x),$ centered at the origin and of radius~$r$
are simply the ellipsoids
$$
{\cal E}_r(0)=\left\{x\in\R^n\colon\quad
\frac{x_1^2}{r^{2\alpha_1}}+\ldots+ \frac{x_n^2}{r^{2\alpha_n}}  <1\right\}
$$
with Lebesgue measure $|\E_r|=C(n) r^\alpha.$  It is easy to see that the unite
sphere with respect to this metric coincides with the unite sphere $\Sigma_n$
with  respect to the Euclidean one.
\begin{defin}\label{d1}
The function $ k(x;\xi)
\colon\R^n\times\{\R^{n}\setminus \{0\}\}\to \R$ is called a {\sl variable
 kernel with mixed homogeneity\/}   if:
\begin{itemize}
\item[$i)$]  for every fixed $x$ the function  $k(x;\cdot)$ is	 a {\sl  constant
kernel\/} satisfying
\begin{itemize}
\item[$i_a)$]
    $k(x;\cdot)\in
C^\infty(\R^{n}\setminus \{0\});$
\item[$i_b)$]
  $k(x;\mu^{\alpha_1}\xi_1,\ldots, \mu^{\alpha_n}\xi_n)=\mu^{-\alpha}k(x;\xi),$
$\forall$ $\mu>0,$ $\alpha_i\geq 1,$ $\alpha=\sum_{i=1}^n\alpha_i;$
\item[$i_c)$]
  $ \ds\int_{\Sigma_n} k(x;\xi)
d\sigma_\xi =0$\quad and\quad $\ds \int_{\Sigma_n}|k(x;\xi)|d\sigma_\xi<\infty;$
\end{itemize}
\item[$ii)$]
 for every multiindex $\beta:$ $\ds \sup\limits_{\xi\in\Sigma_n}
\left|D^\beta_\xi k(x;\xi)\right|\leq C(\beta)$
independently of $x$.
\end{itemize}
\end{defin}
Let us note that in the special case $\alpha_i=1$ and thus $\alpha =n,$
Definition~\ref{d1} gives rise to the classical
Calder\'on-Zygmund kernels.
One more example is when $\alpha_1=\ldots=\alpha_{n-1}=1,$
$\alpha_n=\bar\alpha\geq 1.$  In this case we obtain the kernels studied by
Jones in \cite{BJ} and discussed   in \cite{FR}.

For the sake of completeness we  recall   the definitions and
some properties of the spaces we are going to use.
\begin{defin}\label{dBMO}
For  $f\in L^1_{\rm loc}(\R^n)$ and any ellipsoid $\E\subset\R^n$ centered at
$x\in\R^n$ and of radius $r>0$ set
\begin{equation}\label{modul}
\gamma_f(R):=\sup_{r\leq R} \frac{1}{|{\E}|}\int_{\E}
 |f(y)-f_{\E}|dy\quad  for\ every \ R>0,
 \end{equation}
 where
$f_{\E}=\frac{1}{|\E|}\int_{\E}  f(y) dy $ and
 $|\E|$ is the Lebesgue  measure of
$\E,$ comparable to $r^\alpha.$
Then:
\begin{itemize}
\item[$i)$]  $f\in BMO$ $(${\sl bounded mean oscillation\/}$)$
if
$
\|f\|_\ast:=\sup_R \gamma_f(R)<\infty.$  The quantity
$\|f\|_\ast$ is a norm in $BMO$ modulo constant function under which $BMO$
results a Banach space	$($see \cite{JN}$)$;
\item[$ii)$]  $f\in VMO$ $(${\sl vanishing mean oscillation\/}$)$ with
$VMO$-modulus $\gamma_f(R)$ if
$f$ belongs to $ BMO$  and  $\gamma_f(R)\to 0$ as $R\to 0$ $($see
\cite{S}$)$.
 \end{itemize}
For a bounded domain $\Omega\subset \R^n,$ we	define	   $BMO(\Omega)$
and $VMO(\Omega)$    taking $f\in L^1(\Omega)$ and $\E\cap \Omega$
instead of  $\E$ in $(\ref{modul})$.
\end{defin}
Let $\omega\colon \R^n\times\R_+\to \R_+$ and for any ellipsoid
$\E$
we write $\omega(x,r)=:\omega(\E).$
\begin{defin}\label{dPMS}
A   function $f\in L^p_{\text{loc}}(\R^n),$ $p\in(1,\infty)$ belongs to
the generalized Morrey space $L^{p,\omega}(\R^n)$ if
 the following norm is
finite
\begin{equation}\label{PMS}
\|f\|_{p,\omega}:=\left(\sup_{\E}
\frac{1}{\omega(\E)} \int_{\E}|f(y)|^pdy \right)^{1/p}.
\end{equation}
The space $L^{p,\omega}(\Omega)$ and the norm
$\|f\|_{p,\omega;\Omega}$ are defined by taking $f\in L^p(\Omega)$
and $\E\cap \Omega$ instead of $\E$ in \eqref{PMS}.
\end{defin}
For $\omega(x,r) = 1$ we get  the
Lebesgue
space $L^p(\R^n)$  and for $\omega(x,r)=r^\lambda,$ $\lambda\in(0,\alpha),$
$L^{p,\omega}(\R^n)$ coincides with the Morrey space $L^{p,\lambda}(\R^n)$ when
$\R^n$ is endowed with the  metric $\rho.$
However, there exist weight
functions,  as $\omega(x,r)=r^\lambda\ln(r+2),$ $\lambda\in (0,\alpha)$
 for which $L^{p,\omega}(\R^n)$  does not coincides with any  Morrey
space.

For a given measurable function  $f\in L^1_{\rm loc}(\R^n)$  define  the
{\sl Hardy--Littlewood maximal operator\/} $Mf$  and
the   {\sl sharp maximal operator\/} $f^{\sharp}$   as
$$
Mf(x):=\sup_{x\ni\E}\frac{1}{|\E|}\int_{\E}|f(y)|dy, \quad
f^{\sharp}(x):=\sup_{x\ni\E}\frac{1}{|\E|}\int_{\E}|f(y)-f_{\E}|dy
$$
almost everywhere in $\R^n$ and
 the supremum is taken over all ellipsoids  $\E$  centered at
 $x.$
Define also the operator  $M_sf(x):=(M|f|^s(x))^{1/s}$	for
$s\geq 1.$

The next results are weighted variants of the well-known	maximal and sharp
inequalities  obtained in Lebesgue and Morrey spaces (see  \cite{St}, \cite{CF},  \cite{DPR}).
\begin{lem}\label{l1} {\sc  (Maximal inequality)(\cite{Na})}
  Assume that there are
constants $C_1$  and $C_2$  such that for any $x_0\in \R^n$  and for any $r>0$
\begin{align}\label{Na1}
r\leq t\leq 2r \Longrightarrow	   C_1\leq \frac{\omega(x_0,t)}{\omega(x_0,r)}
\leq C_2,\\
\label{Na2}
  \int_r^\infty   \frac{\omega(x_0,t)}{t^{\alpha+1}} dt\leq
C\frac{\omega(x_0,r)}{r^\alpha}.
\end{align}
 For $1\leq s<p<\infty,$  there is a constant $C_{p,s}>0 $  such that for
$f\in L^{p,\omega}(\R^n)$
 $$
\|M_sf\|_{p,\omega}\leq  C_{p,s} \|f\|_{p,\omega}.
$$
\end{lem}
\begin{lem}\label{l2} {\sc  (Sharp inequality)}
Let $0<\sigma\leq 1$  and $\E$	be an ellipsoid centered at $x_0\in \R^n$ of
radius $r.$ Suppose that $\omega(x_0,r)$ satisfies \eqref{Na1}
and
$$
\int_r^\infty \frac{\omega(x_0,t)}{t^{\sigma \alpha+1}} dt\leq
C\frac{\omega(x_0,r)}{r^{\sigma\alpha}}.
$$
Then for $p\in(1,\infty)$
and  $f\in L^{p,\omega}(\R^n)$ there exists a constant $C$ independent of $f$
such that
$$
\|f\|_{p,\omega}\leq  C \|f^{\sharp}\|_{p,\omega}.
$$
\end{lem}
\begin{pf}
Let $\chi_{\E}$  be the characteristic function of the ellipsoid and denote
by $2\E$    an ellipsoid centered  at $x_0$  and of radius $2r.$ It is easy to
verify that
$
M\chi_{\E}(x) \leq r^\alpha/(\rho(x-x_0)-r)^\alpha \leq 1,$ for all $x\in\R^n. $
 Further, for any $x\in 2^{k+1}\E\setminus 2^k\E,  $ $k=1,2,\ldots$  the maximal
function of $\chi_{\E}$    could be estimated by
$
r^\alpha / (2^{k+1}r-r)^\alpha	\leq M\chi_{\E}(x) \leq
r^\alpha / (2^{k}r-r)^\alpha
$
which gives a reason to compare        $M\chi_{\E}(x)$ with $2^{-k\alpha}$  for
any $x$  as above.
From the properties
 of the maximal function, that is $|f|\leq Mf$	and $Mf\leq
f^\#$ (see \cite[p. 410]{GR})	follows

\begin{align*}
J&=\int_{\E}|f(y)|^p
dy=\int_{\R^{n}}|f(y)|^p \chi_{\E}(y)dy\\
&\leq\int_{\R^{n}}|Mf(y)|^p
\big(M\chi_{\E}(y)\big)^{\sigma}dy
\leq C\int_{\R^{n}}|f^{\#}(y)|^p
\big(M\chi_{\E}(y)\big)^{\sigma}dy\\
&\leq C\Big\{\int_{2\E}|f^{\#}(y)|^pdy\\
&\qquad +\sum_{k=1}^{\infty}\int_{2^{k+1}\E \setminus 2^k\E  }
|f^{\#}(y)|^p\Big(\frac{r}{\rho(y-x_0)-r}
\Big)^{\sigma \alpha}dy\Big\}\\
&\leq C\Big\{ \omega(2\E) \frac{1}{\omega(2\E)}
\int_{2\E}|f^{\#}(y)|^pdy\\
&\qquad + r^{\sigma \alpha}\sum_{k=1}^\infty \frac{
\omega(2^{k+1}\E)}{(2^{k}r)^{\sigma \alpha}}
\frac{1}{\omega(2^{k+1}\E)} \int_{2^{k+1}\E} |f^{\#}(y)|^p dy
 \Big\}\\
& \leq C r^{\sigma \alpha} \sum_{k=0}^\infty
\frac{\omega(2^k\E)}{(2^{k}r)^{\sigma \alpha}} \|f^\#\|_{p,\omega}^p.
 \end{align*}
 From the properties of the function $\omega(x_0,t)$ follows
$$
\frac{\omega(2^k\E)}{(2^{k}r)^{\sigma \alpha}} \sim
\int_{2^k r}^{2^{k+1}r} \frac{\omega(x_0,t)}{t^{\sigma
\alpha+1}}dt.
$$
Hence
$$
\int_\E |f(y)| dy\leq C r^{\sigma \alpha} \int_r^\infty
\frac{\omega(x_0,t)}{t^{\sigma \alpha+1}} dt \|f\|^p_{p,\omega} \leq
C\omega(x_0,r) \|f\|^p_{p,\omega}.
$$
\end{pf}
\begin{lem}\label{l3} {\sc (John-Nirenberg type lemma)}
Let   $a\in BMO$ and  $p\in(1,\infty).$ Then
for any ellipsoid $\E$	holds
$$
\left(\frac{1}{|\E|} \int_{\E}\left|a(y)-a_{\E}\right|^pdy
\right)^{1/p}\leq C(p)\|a\|_\ast.
$$
\end{lem}

One more background we need is that for {\sl spherical harmonics\/}  and their
properties (see for instance  \cite{CZ}, \cite{FR}, \cite{CFL}).
Recall that any homogeneous polynomial $P\colon\ \R^n\to\R$
of  degree $m$ that satisfies  $\Delta P(x)=0$ is called an {\sl $n$-dimensional
solid harmonic of
degree $m.$\/}	 Its restriction to the unit sphere $\Sigma_n$	will be
called an  {\sl $n$-dimensional  spherical harmonic of degree $m.$\/}
Denote by $\Upsilon_m$ the space of all $n$-dimensional spherical  harmonics
of
degree $m.$
 In general it	results a finite-dimensional linear
space with  $g_m=\dim\Upsilon_m$ such that  $g_0=1,$ $g_1=n$   and
\begin{equation}\label{*}
g_m={m+n-1\choose n-1}-{m+n-3 \choose n-1} \leq C(n) m^{n-2}, \   m\geq 2.
\end{equation}
 Further, let $\{Y_{sm}(x)\}_{s=1}^{g_m}$ be an
{\sl orthonormal base\/} of $\Upsilon_m.$ Then
 $\{Y_{sm}\}_{s=1}^{g_m}{}_{m=0}^{\infty}$ is a {\sl complete orthonormal
system\/} in $L^2(\Sigma_n)$ and
\begin{equation}\label{Y1}
\sup\limits_{x\in\Sigma_n} \left|D^\beta_x  Y_{sm}(x)  \right| \leq C(n)
m^{|\beta|+(n-2)/2},\quad m=1,2,\ldots.
\end{equation}
If, for instance,  $\phi\in C^\infty(\Sigma_n)$
then $\sum_{s,m} b_{sm} Y_{sm}(x) $  is the Fourier series
expansion of $\phi(x)$ with respect to $\{Y_{sm}(x)\}_{s,m}$
$\left(\sum_{s,m}\right.$
substitutes $\left.\sum_{m=0}^{\infty} \sum_{s=1}^{g_m}\right)$ and
\begin{equation}\label{Y3}
b_{sm}= \int_{\Sigma_n}\phi ( y) Y_{sm}( y)d\sigma,\quad
|b_{sm}|\leq C(n,l)
m^{-2l}\sup\limits_{\overset{|\beta|=2l}{y\in\Sigma_n}}
\left|D^\beta_y\phi(y) \right|
\end{equation}
for any integer $l.$ In particular, the expansion of $\phi$ into spherical
harmonics converges uniformly to $\phi.$ For the proof of the above results see
\cite{CZ}.

\section{Singular integral estimates}\label{s3}

\setcounter{equation}{0}
\setcounter{thm}{0}

Let $k(x;\xi)$	be a	kernel in the sense of Definition~\ref{d1}.
In order to ensure the existence of the operators  \eqref{Kf} and \eqref{Cakf}	in
$L^p(\R^n)$ we	restrict our considerations to functions $f\in L^p(\R^n),$
$1<p<\infty$  for which the norm  \eqref{PMS}  is finite. For the sake of
convenience we still denote these spaces by $L^{p,\omega}(\R^n).$
Having	 in mind this  we define the operators
$\K_\varepsilon f$ and	$\C_\varepsilon [a,k] f$ for   $a\in BMO$ and
$f\in  L^{p,\omega}(\R^n)$ with $p\in(1,\infty)$ and $\omega$ satisfying
\eqref{Na1} and \eqref{Na2},
 by
\begin{align*}
\K_\varepsilon f(x):=&\ \ds \int_{\rho(x-y)>\varepsilon}
 k(x;x-y)f(y)dy,\\
\C_\varepsilon [a,k]f(x):=&\ \ds \K_\varepsilon (af)(x) -a(x)\K_\varepsilon
f(x)\\
\phantom{:}=&\	 \ds  \int_{\rho(x-y)>\varepsilon}
	   k(x;x-y)[a(y)-a(x)]f(y)dy.
\end{align*}
We are going to prove  that $\K_\varepsilon$ and $\C_\varepsilon [a,k]$
are bounded and continuous  from $ L^{p,\omega}(\R^n)$ into itself  uniformly in
$\varepsilon.$
This along with the properties of the kernel $k(x;\xi)$ will enable to let $\varepsilon
\to0$ obtaining as limits in the $L^{p,\omega}(\R^n)$-topology the {\sl
singular\/} integrals
\begin{align*}
\K f(x):=&\ \ds P.V. \int_{\R^n}
k(x;x-y)f(y)dy=\lim_{\varepsilon \to 0} \K_\varepsilon f(x)\\
\C[a,k]f(x):=&\ \ds P.V.\int_{\R^n} k(x;x-y)[a(y)-a(x)]f(y)dy=
\lim_{\varepsilon \to 0} \C_\varepsilon[a,k]f(x).
\end{align*}
Moreover, we shall show that the last ones  are also continuous in
$L^{p,\omega}(\R^n).$

Let us note  assuming $f\in L^p(\R^n),$ $p\in(1,\infty)$ Fabes-Rivi\`ere
(\cite{FR}) show that $\K f$ exists in $L^p(\R^n)$ for $p\in(1,\infty)$ as a
limit of $\K_\varepsilon f$ when $\varepsilon\to 0$ in the $L^p$-norm.
Moreover, the operator
$\K\colon\ L^p(\R^n)\to L^p(\R^n)$ is continuous and this leads also to
continuity in $L^p(\R^n)$ of $\C[a,k]f$ if $a(x)$ is essentially bounded.
As it concerns to the commutator we are going to derive a result similar
to
that of Coifman-Rochberg-Weiss (\cite[Theorem~1]{CRW}), which asserts: if $\K$ is
Calder\'on-Zygmund operator in $L^p(\R^n),$ $p\in(1,\infty)$   and $a\in BMO$
than the commutator $\C[a,\cdot]$  is a well defined linear continuous operator
from $L^p(\R^n)$   into itself.   Later, this result has been extended	 by
Bramanti-Cerutti (\cite{BC}) in the framework of homogeneous spaces.  Based on
this background about Calder\'on-Zygmund operators, we are going
 to obtain
continuity in $L^{p,\omega}(\R^n)$   and boundedness in terms of $\|a\|_\ast$
for the commutators  \eqref{Cakf}   having  kernel of more general type.

\begin{thm}\label{th1}	Let $k(x;\xi)$	be a variable kernel of mixed
homogeneity, $f\in  L^{p,\omega}(\R^n),$ $p\in (1,\infty),$ $\omega$ satisfies
\eqref{Na1} and \eqref{Na2}, and $a\in	BMO.$
Then there exist the
 integrals $\K f,\
\C[a,k]f\in L^{p,\omega}(\R^n)$
as limits of
$\K_\varepsilon  f$ and $\C_\varepsilon [a,k]f$ when $\varepsilon
\to0$ with respect to the $L^{p,\omega}(\R^n)$-norm.
 The operators $\K$ and $\C [a,k]$ are bounded from
$ L^{p,\omega}(\R^n)$ into itself and
$$
\|\K f\|_{p,\omega}   \leq C\|f\|_{p,\omega},\qquad
\|\C[a,k]f \|_{p,\omega} \leq C\|a\|_\ast \|f\|_{p,\omega}
$$
where the constants depend on $n,p,\alpha$  and  $k$ through the
constant $C(\beta).$
\end{thm}
\begin{pf}
Let $x,y\in \R^{n}$  and $\ol{y}={y}/{\rho(y)}\in \Sigma_n.$  From
the properties of the kernel with respect to the second variable and
  the completeness of
$\{Y_{sm}(x)\}_{s,m}$  in
    $L^2(\Sigma_n  )$ it follows
$$
k(x;x-y)=\rho(x-y)^{-\alpha} k(x;\ol{x-y})=\rho(x-y)^{-\alpha}\sum_{s,m}
b_{sm}(x) Y_{sm}(\ol{x-y}).
$$
This way, the  Definition~\ref{d1}~$ii)$ and
\eqref{Y3} imply
\begin{equation}\label{Y4}
\|b_{sm}\|_\infty  \leq C(n,l,k)m^{-2l}
\end{equation}
for any integer $l>1.$
Replacing the kernel with its expansion, we get
\begin{align}\label{K_eps}
\K_\varepsilon f(x)=&\ \ds
\int_{\rho(x-y)>\varepsilon}
\sum_{s,m}b_{sm}(x)
\H_{sm}(x-y) f(y)dy,\\
\nonumber
\C_\varepsilon [a,k]f(x) =&\ \ds
\int_{\rho(x-y)>\varepsilon}
\sum_{s,m}b_{sm}(x)
\H_{sm}(x-y) [a(y) - a(x)]f(y)dy
\end{align}
with $\H_{sm}(x-y)$ standing for   $Y_{sm}(\ol{x-y})
\rho(x-y)^{-\alpha}.$ It is easy to check that $\H_{sm}(\cdot)$  is a constant
kernel in the sense of Definition~\ref{d1}~$i)$. Indeed, $i_a)$ and $i_b)$ are
trivial while $i_c)$ follows from the fact that $Y_{sm}(x)$ is a harmonic
homogeneous polynomial	and the property of integral mean on
sphere for  harmonic functions (i.e.  $Y_{sm}(0)=0$).  In order to get series expansions
of $\K_\varepsilon f$ and $\C_\varepsilon[a,k]f,$ we let $x\in \R^n$ and
$y\in\R^n$ to be such that $\rho(x-y)>\varepsilon.$ Then \eqref{*}, \eqref{Y1}
and \eqref{Y4} yield
$$
\left|\sum_{m=1}^N \sum_{s=1}^{g_m} b_{sm}(x)\frac{Y_{sm}(\overline{x-y})}{
\rho(x-y)^\alpha} f(y)\right| \leq C(n) \frac{|f(y)|}{\rho(x-y)^\alpha}
\sum_{m=1}^{\infty} m^{n-2+(n-2)/2-2l}
$$
where $|f(\cdot)|\rho(x-\cdot)^{-\alpha}\in L^1(\R^n)$ for a.a. $x\in \R^n$
and  the integer $l$ is preliminary chosen greater
than
$(3n-4)/4.$ Similar inequality holds also for
the commutator	$\C_\varepsilon[a,k]f.$ Thus,
by the dominated convergence theorem
\begin{align}\label{EK}
\K_\varepsilon f(x)= & \sum_{s,m}b_{sm}(x) \K_{sm,\varepsilon}(x),\\
\nonumber
\C_\varepsilon [a,k]f(x) = &
     \sum_{s,m}b_{sm}(x) \C_{sm,\varepsilon}[a,k]f(x)
\end{align}
 with
\begin{align*}
\K_{sm,\varepsilon}(x):= & \int_{\rho(x-y)>\varepsilon}
	   \H_{sm}(x-y) f(y)dy,\\
\C_{sm,\varepsilon} [a,k]f(x):= & \int_{\rho(x-y)>\varepsilon}
		    \H_{sm}(x-y) [a(y) - a(x)]f(y)dy.
 \end{align*}
This way instead of the operators $\K f$  and $\C[a,k]f$ we shall study the
existence and boundedness in $L^{p,\omega}(\R^n)$  of the singular integrals
\begin{align*}
\K_{sm}f(x):= & P.V. \ds\int_{\R^{n}}\H_{sm}(x-y)  f(y)dy,\\
 \C_{sm}[a,k]f(x):= & P.V \ds\int_{\R^{n}}
\H_{sm}(x-y)[a(y)-a(x)] f(y) dy
\end{align*}
with constant kernels $\H_{sm}(\cdot).$
For what concern boundedness of $\K_{sm}$ in $L^p(\R^n)$ we dispose of
\cite[Theorem~II.1]{FR} and this implies, through \cite[Theorem~2.5]{BC},
boundedness in $L^p(\R^n)$ of $\C_{sm}[a,k]$ as well.
 The cited results however require the kernel to
have some ``integral continuity'',
 called {\sl the H\"ormander  condition\/}. It turns out that
$\H_{sm}(\cdot)$ satisfies even stronger condition as shows the following lemma.
\begin{lem}\label{l5} {\sc (Pointwise H\"ormander condition)}
Let $\E$ and $2\E$ be  ellipsoids  centered at	$x_0$  and of radius $r$
and $2r,$ respectively. Then
\begin{equation}\label{PHE}
\left|\H_{sm}(x-y)- \H_{sm}(x_0-y)  \right|
 \leq C(n,\alpha) m^{n/2}
\frac{\rho(x_0-x)}{\rho(x_0-y)^{\alpha+1}}
\end{equation}
for each $x\in \E$  and $y\notin 2\E.$
\end{lem}
\begin{pf}
We shall apply the mean value theorem to $\H_{sm}$ and therefore decay estimate
for the gradient $\nabla \H_{sm}$ is needed. Let $x\in\R^n\setminus\{0\}$
be
an arbitrary point. The implicit function theorem applied to the equation
$F(x,\rho(x))=1$ gives an  expression for  the gradient $\nabla \rho(x)$
and
straightforward calculations imply
\begin{align*}
\ds \frac{\partial\H_{sm}}{\partial x_i}(x) =
\ds \frac{1}{\rho(x)^{\alpha+\alpha_i}}&\Big( -\alpha Y_{sm}(\ol{x})
\frac{x_i}{\rho(x)^{\alpha_i}
  \sum_{j=1}^n \alpha_j x_j^2 \rho(x)^{-2\alpha_j}}
\ds + \frac{\partial Y_{sm}}{\partial \ol{x}_i}(\ol{x})\\
\ds & - \sum_{k=1}^n \alpha_k \frac{\partial
\ds Y_{sm}}{\partial  \ol{x}_k}(\ol{x})
\frac{x_ix_k}{\rho(x)^{\alpha_i}  \rho(x)^{\alpha_k}
 \sum_{j=1}^n \alpha_j x_j^2 \rho(x)^{-2\alpha_j}}\Big).
\end{align*}
Since $\ol{x}\in\Sigma_n$ and taking
into account \eqref{Y1}, $x_i /\rho(x)^{\alpha_i}\leq |\ol x|\leq 1$ and
$\min\alpha_i \leq \sum_{j=1}^n \alpha_j x_j^2 / \rho(x)^{2\alpha_j}\leq
\max\alpha_i,$ we get
\begin{equation}\label{Hx}
\left|\frac{\partial\H_{sm}}{\partial x_i}(x)\right|
\leq C(n,\alpha) \frac{m^{n/2}}{\rho(x)^{\alpha+\alpha_i}}\qquad \forall\
x\in \R^n\setminus\{0\}.
\end{equation}
Now, applying the mean value theorem to the left-hand side of \eqref{PHE}  we
get
\begin{equation}\label{gradH}
\H_{sm}(x-y)- \H_{sm}(x_0-y) =\sum_{i=1}^n \frac{\partial \H_{sm}}{\partial
x_i} (x_0-\xi) (x_0-x)_i
\end{equation}
with $\xi=y-t(x-x_0)$ and $t\in(0,1).$	Obviously
$\rho(y-\xi)=t\rho(x_0-x)\leq r$  which along with $y\not\in 2\E$   gives  that
$\xi$  does not belong to $\E$	and $\rho(x_0-\xi)\geq \frac{1}{2}\rho(x_0-y).$
 Having in mind $(x_0-x)_i\leq
\rho(x_0-x)^{\alpha_i},$ \eqref{Hx} and
 \eqref{gradH} we
obtain
\begin{align*}
&\left|\H_{sm}(x-y)- \H_{sm}(x_0-y)  \right| \leq C(n,\alpha) m^{n/2}
\sum_{i=1}^n \frac{\rho(x_0-x)^{\alpha_i}}{\rho(x_0-\xi)^{\alpha+\alpha_i}}\\
&\quad \leq  C(n,\alpha) m^{n/2} \frac{\rho(x_0-x)}{\rho(x_0-y)^{\alpha+1}}
\sum_{i=1}^n \frac{\rho(x_0-x)^{\alpha_i-1}}{\rho(x_0-\xi)^{\alpha_i-1}}\\
&\quad \leq C(n,\alpha) m^{n/2}
\frac{\rho(x_0-x)}{\rho(x_0-y)^{\alpha+1}}
\end{align*}
where we have used that $\alpha_i\geq 1$ and
$\rho(x_0-x)<\frac{1}{2}\rho(x_0-y)\leq  \rho(x_0-\xi)$  from which follows
immediately  the last sum is no greater than $n.$
\end{pf}
\begin{rem}\label{remarkHIC}\em
This result ensures
   the kernel $\H_{sm}$ satisfies
the H\"ormander integral condition (see  \cite[(1.1)]{FR})
$$
\int_{\{y\in\R^n\colon \rho(y)\geq 4\rho(x)\}}
\left|\H_{sm}(y-x)-\H_{sm}(y)\right|dy \leq C
$$
with a constant  independent of $x.$
\end{rem}

In view of the cited above results  there exist $\K_{sm}f,\ \C_{sm}[a,k]f\in L^p(\R^n)$
such
that $$
\lim_{\varepsilon\to0}
\left\| \K_{sm,\varepsilon}f-\K_{sm}f\right\|_{L^p(\R^n)}=
\lim_{\varepsilon\to0}
\left\| \C_{sm,\varepsilon}[a,k]f-\C_{sm}[a,k]f\right\|_{L^p(\R^n)}=0.
$$
Our goal is to show that this convergence is fulfilled also with respect to the
$L^{p,\omega}(\R^n)$-norm. The proof is broken up into several Lemmas.
\begin{lem}\label{l6} The singular integrals
$\K_{sm}f$ and their commutators $\C_{sm}[a,k]f$ satisfy
 \begin{align}\label{Ksharp}
 (\K_{sm}f)^{\sharp}(x)\leq &	C m^{n/2}
\big(M(|f|^p)(x)\big)^{1/p},\\
\nonumber
 (\C_{sm}[a,k]f)^{\sharp}(x)\leq &   C
\|a\|_\ast\Big\{\big(M(|\K_{sm}f|^p)(x)\big)^{1/p}
 +  m^{n/{2}}
\big(M(|f|^p)(x)\big)^{1/p} \Big\},
\end{align}
where the constant depends on $n,$ $p$ and   $\alpha$ but not
on  $f.$
\end{lem}
\begin{pf}
For arbitrary $x_0\in \R^{n},$ set  $\E$ for the ellipsoid $\E$
centered at $x_0$ and of radius $r.$
Let we consider the  expression
$$
I:=\frac{1}{|\E|} \int_{\E}
\left|\K_{sm}f(y)  -(\K_{sm}f)_{\E}\right|dy.
$$
Adding and extracting $\K_{sm,2r}f(x_0)$
to the function under the sign of the integral we obtain
$$
I \leq \frac{2}{|\E|} \int_{\E}
\left|\K_{sm}f(y) -\K_{sm,2r}f(x_0)\right| dy :=2 I(x_0,\E).
$$
Set $(2\E)^c=\R^{n}\setminus 2\E$ and write
$f=f\chi_{2\E}+f\chi_{(2\E)^c}=f_1+f_2$ with $\chi$ being the characteristic
function of the respective set.  Hence
\begin{align*}
I(x_0,\E)\leq&\ \frac{1}{|\E|} \int_{\E}
|\K_{sm}f_1(y)| dy\\
 &\ + \frac{1}{|\E|} \int_{\E}
|\K_{sm}f_2(y) -\K_{sm,2r}f(x_0)|  dy=:I_1(x_0,\E)+I_2(x_0,\E).
\end{align*}
From the boundedness of
$\K_{sm}$ in $L^p(\R^n)$
(\cite[Theorem~II.1]{FR}) follows
\begin{align*}
I_1(x_0,\E)\leq&\ \frac{1}{|\E|}\left(\int_{\E}1 dy
\right)^{1/p'}\left(\int_{\E}
|\K_{sm}f_1(y)|^pdy \right)^{1/p}=\frac{1}{|\E|^{1/p}} \|\K_{sm}f_1\|_p\\
\leq&\ \frac{C(p,\alpha)}{|\E|^{1/p}} \|f_1\|_p  \leq C(p,\alpha)
\big(M(|f|^p)(x_0)\big)^{1/p}
\end{align*}
with $1/p'+1/p=1.$
About $I_2(x_0,\E),$ we have

 \begin{align*}
I_2(x_0,\E)\leq&\ \frac{1}{|\E|} \int_{\E}\left( \int_{(2\E)^c} \left|
\H_{sm}(y-\xi) -\H_{sm}(x_0-\xi)     \right| |f(\xi)| d\xi\right) dy\\
\leq&\ C(n,\alpha)m^{n/2} \frac{1}{|\E|} \int_{\E}
\left(\ \int_{(2\E)^c} \frac{\rho(x_0-y)}{\rho(x_0-\xi)^{\alpha+1}}
|f(\xi)| d\xi\right) dy\\
\leq&\	C(n,\alpha)m^{n/2} r \sum_{k=1}^{\infty}
\int_{2^{k+1}\E\setminus 2^k\E} \frac{|f(\xi)|}{\rho(x_0-\xi)^{\alpha+1}}
d\xi\\
\leq&\	C(n,\alpha)  m^{n/2} \frac{1}{r^{\alpha}} \sum_{k=1}^{\infty}
\frac{1}{2^{k(\alpha+1)}}|2^{k+1}\E|
 \left(\frac{1}{|2^{k+1}\E|}
\int_{2^{k+1}\E} |f(\xi)|^pd\xi \right)^{1/p}\\
\leq&\	C(n,\alpha) m^{n/2} \big(M(|f|^p)(x_0)\big)^{1/p},
\end{align*}
after applying Lemma~\ref{l5} for  $y\in\E$ and $\xi\in (2\E)^c.$
Taking $\sup_{\E}I(x_0,\E)$ and heaving in mind the arbitrarity
of $x_0,$ we obtain \eqref{Ksharp} for any $x\in \R^{n}.$

To estimate the sharp function of	the commutator we shall employ	the
idea of Stromberg (see \cite{To}) which consists of expressing
$\C_{sm}[a,k]f$
as a sum of  integral operators and estimating	their sharp functions.
Precisely,
\begin{align*}
&\C_{sm}[a,k]f(x)= \K_{sm}(a   -a_{\E})f(x)-(a(x)-a_{\E})\K_{sm}f(x)\\
&\qquad = \K_{sm}(a	-a_{\E})f_1(x) + \K_{sm}(a   -a_{\E})f_2(x)
-(a(x)-a_{\E})\K_{sm}f(x) \\
&\qquad =: J_1(x)+J_2(x) +J_3(x)
\end{align*}
where we have used the same truncation for the function $f$ as in $I(x_0,\E).$
Before proceed further, let us point out the obvious inequality
\begin{equation}\label{a1}
|a_{2\E}-a_{\E}|\leq C(n,\alpha)\|a\|_\ast\qquad \forall\ a\in BMO(\R^n)
\end{equation}
and its by-product
\begin{equation}\label{a2}
|a_{2^k\E}-a_{\E}|\leq C(n,\alpha)k\|a\|_\ast
\end{equation}
following from \eqref{a1} by running induction.
Now, for arbitrary $p\in(1,\infty)$  and $q\in (1,p),$	we have

\begin{align*}
G_1(x_0,\E)&:= \frac{1}{|\E|}\int_{\E} |J_1(x)-(J_1)_{\E}| dx
\leq\frac{2}{|\E|} \int_\E |\K_{sm}(a	-a_{\E})f_1(x)|dx\\
 \leq&\ \frac{2}{|\E|}\left(\	\int_{\E}
|\K_{sm}(a   -a_{\E})f_1(x)|^qdx  \right)^{1/q} \left(\  \int_{\E} 1
dx \right)^{1/q'}\\
\leq&\ \frac{C(q,\alpha)}{|\E|^{1/q}} \left(\ \int_\E
\big|(a(x)-a_\E)f_1(x)\big|^q dx \right)^{1/q}\\
 \leq&\ \frac{C(q,\alpha)}{|\E|^{1/q}} \left(\	\int_{2\E}|f(x)|^pdx
\right)^{1/p} \left(\  \int_{2\E}|a(x)-a_{\E}|^{pq/(p-q)} dx\right)^{(p-q)/pq}.
\end{align*}
 Further, \eqref{a1} and Lemma~\ref{l3} applied to the second integral yield
\begin{align*}
&\int_{2\E} |a(x)-a_\E|^{pq/(p-q)}dx
   \leq C(p,q)\left(\
       \int_{2\E} |a(x)-a_{2\E}|^{pq/(p-q)}dx\right.\\
&\qquad\qquad \qquad\qquad\qquad
\left.+  \int_{2\E} |a_{2\E}-a_\E|^{pq/(p-q)}dx\right)\\
&\quad\leq C(p,q)\left(|2\E|\frac{1}{|2\E|}
       \int_{2\E} |a(x)-a_{2\E}|^{pq/(p-q)}dx+
       |2\E| C(n,\alpha) \|a\|_{\ast}^{pq/(p-q)}\right)\\
&\quad\leq C(n,p,q,\alpha)|2\E|\|a\|_{\ast}^{pq/(p-q)}.
 \end{align*}
Therefore,
$$
G_1(x_0,\E)\leq C \|a\|_\ast \left(\frac{1}{|2\E|}\int_{2\E}|f(y)|^pdy
\right)^{1/p}\leq C\|a\|_\ast \big(M(|f|^p)(x_0)\big)^{1/p}.
 $$
To estimate the sharp function of $J_2(x),$ we proceed analogously as we already
did for  $I_2(x_0,\E).$ Precisely,
$$
G_2(x_0,\E):=\frac{1}{|\E|} \int_{\E} |J_2(x)-(J_2)_{\E}|dx \leq
 \frac{2}{|\E|}\int_{\E} |J_2(x)-J_2(x_0)| dx
$$
and the integrand satisfies
\begin{align*}
&|J_2(x)-J_2(x_0)| \leq  \int_{(2\E)^c}\left|  \H_{sm}(x-y)-
\H_{sm}(x_0-y)\right|  |a(y)-a_{\E}|\, |f(y)| dy\\
&\ \leq C(n,\alpha) m^{n/2}\rho(x_0-x) \int_{(2\E)^c}
\frac{|a(y)-a_{\E}| |f(y)|}{\rho(x_0-y)^{\alpha+1}}dy\\
&\  \leq C(n,\alpha) m^{n/2} r \left(\ \int_{(2\E)^c}
\frac{|f(y)|^p}{\rho(x_0-y)^{\alpha+1}}dy\right)^{1/p}
 \left(\ \int_{(2\E)^c}
\frac{|a(y)-a_{\E}|^{p'}}{\rho(x_0-y)^{\alpha+1}}dy\right)^{1/p'}
\end{align*}
where $1/p+1/p'=1$ and we have applied the H\"ormander
pointwise estimate (Lemma~\ref{l5}).
Later on,
$$
\int_{(2\E)^c}
\frac{|f(y)|^p}{\rho(x_0-y)^{\alpha+1}}dy  =
\sum_{k=1}^{\infty}
\int_{2^{k+1}\E\setminus 2^k\E}
\frac{|f(y)|^p}{\rho(x_0-y)^{\alpha+1}}dy
\leq\frac{2^{\alpha+1}}{r}M(|f|^p)(x_0),
$$
while \eqref{a2} and Lemma~\ref{l3} imply
\begin{align*}
&\int_{(2\E)^c}
\frac{|a(y)-a_{\E}|^{p'}}{\rho(x_0-y)^{\alpha+1}}dy= \sum_{k=1}^{\infty}
\int_{2^{k+1}\E\setminus 2^k \E}
\frac{|a(y)-a_{\E}|^{p'}}{\rho(x_0-y)^{\alpha+1}}dy\\
&\quad\leq \sum_{k=1}^{\infty} \frac{1}{(2^k r)^{\alpha+1}}
\int_{2^{k+1}\E}|a(y)-a_{\E}|^{p'}dy\\
&\quad \leq \sum_{k=1}^\infty \frac{2^{p'-1}}{(2^kr)^{\alpha+1}}
\int_{2^{k+1}\E}\left(|a(y)-a_{2^{k+1}\E}|^{p'}+ |a_{2^{k+1}\E}-a_\E|^{p'}
\right) dy\\
&\quad \leq C\sum_{k=1}^\infty \frac{|2^{k+1}\E|}{(2^kr)^{\alpha+1}}
(1+k^{p'})\|a\|_\ast^{p'}
\leq  C \frac{\|a\|_\ast^{p'}}{r}
\end{align*}
and the constant depends on $n,$ $p$ and $\alpha.$
Hence
$$
G_2(x_0,\E)\leq C(n,p,\alpha) m^{n/2}	\|a\|_\ast \big(M(|f|^p)(x_0)\big)^{1/p}.
$$
Finally,
\begin{align*}
G_3(x_0,\E):=&\ \frac{1}{|\E|} \int_{\E} |J_3(x)-(J_3)_{\E}| dx
\leq   \frac{2}{|\E|} \int_{\E} |a(x)-a_{\E}| |\K_{sm}f(x)|dx\\
\leq&\ 2\left(\ \frac{1}{|\E|}\int_{\E} |a(x)-a_{\E}|^{p'}dx\right)^{1/p'}
\left(\  \frac{1}{|\E|}\int_{\E} |\K_{sm}f(x)|^{p}dx\right)^{1/p}\\
\leq&\ C(p) \|a\|_\ast \big(M(|\K_{sm}f|^p)(x_0)\big)^{1/p}.
\end{align*}
Summing up $G_1(x_0,\E),$  $G_2(x_0,\E)$ and  $G_3(x_0,\E)$ and taking
 the supremum with respect to
$\E$ and rendering in account the arbitrarity of the point $x_0$  we get the desired  estimate for the commutator.
\end{pf}
\begin{lem}\label{l6a}
The operators  $\K_{sm}$  and  $\C_{sm}[a,k]$ are continuous acting
 from $ L^{p,\omega}(\R^n)$ into itself and
\begin{equation}\label{Kpl}
\|\K_{sm}f\|_{p,\omega} \leq C m^{n/2}\|f\|_{p,\omega},\
\|\C_{sm}[a,k]f\|_{p,\omega} \leq C
m^{n/2}\|a\|_\ast\|f\|_{p,\omega}
\end{equation}
with constants depending on $n,$  $p,$	and $\alpha.$
\end{lem}
\begin{pf} First of all we shall   estimate  the $L^{p,\omega}$-norms
of the	sharp functions of the considered  operators. Since
 the expression for $(\K_{sm}f)^\sharp$ in
\eqref{Ksharp} holds true for
any $q\in(1,p)$ as well, the maximal inequality (Lemma~\ref{l1}) with
$s=1$ asserts
\begin{align*}
\int_\E |(\K_{sm}f)^{\sharp}(x)|^p dx & \leq Cm^{pn/2}
  \int_\E |M(|f|^q)(x)|^{p/q}dx\\
&  \leq Cm^{pn/2} \omega(\E) \|M(|f|^q)\|^{p/q}_{p/q,\omega}\\
&\leq Cm^{pn/2} \omega(\E)   \| |f|^q \|^{p/q}_{p/q,\omega}
 \leq Cm^{pn/2} \omega(\E)   \| f \|^{p}_{p,\omega}.
\end{align*}
Dividing of $\omega(\E)$ and taking $\sup_{\E},$ we arrive at
$$
\|(\K_{sm}f)^{\sharp}  \|_{p,\omega}\leq C m^{n/2}
\|f\|_{p,\omega}
$$
which implies the first inequality in  \eqref{Kpl} through Lemma~\ref{l2}.
The  $L^{p,\omega}$-estimate  for the commutator follows in the
same
manner.
\end{pf}
\begin{lem}\label{c1}  The  operators
$\K_{sm,\varepsilon}$  and   $\C_{sm,\varepsilon}[a,k]$
 are continuous  acting from $L^{p,\omega}(\R^n)$ into itself and satisfy
\begin{equation}\label{Kpl1}
\|\K_{sm,\varepsilon}f\|_{p,\omega} \leq C
m^{n/2}\|f\|_{p,\omega},\
\|\C_{sm,\varepsilon}[a,k]f\|_{p,\omega} \leq C
m^{n/2}\|a\|_\ast\|f\|_{p,\omega}
\end{equation}
with constants depending on $n,$ $p$ and
$\alpha.$
\end{lem}
\begin{pf}
Let  $\E_\varepsilon$ and $\E_{\varepsilon/2}$ be ellipsoids  centered at
$x\in\R^n$ and of radius $\varepsilon$ and $\varepsilon/2,$ respectively.
Writing $f=f\chi_{\E_\varepsilon}+f\chi_{(\E_\varepsilon)^c}=f_1+f_2$  we obtain
\begin{align*}
&\K_{sm,\varepsilon} f(x)\leq
\frac{1}{|\E_{\varepsilon/2}|}\int_{\E_{\varepsilon/2}}  |\K_{sm,\varepsilon}f(x)|dy
\leq  \frac{1}{|\E_{\varepsilon/2}|}
\int_{\E_{\varepsilon/2}}|\K_{sm}f(y)|dy\\
&\qquad  +\frac{1}{|\E_{\varepsilon/2}|}\int_{\E_{\varepsilon/2}}
|\K_{sm,\varepsilon}f(x)-\K_{sm} f(y)|dy\\
 &\leq \frac{2}{|\E_{\varepsilon/2}|}\int_{\E_{\varepsilon/2}}
|\K_{sm}f_1(y)|dy + \frac{1}{|\E_{\varepsilon/2}|}\int_{\E_{\varepsilon/2}}
|\K_{sm,\varepsilon}f(x)-\K_{sm} f_2(y)| dy\\
& := 2I_1(x,\E_{\varepsilon/2})+I_2(x,\E_{\varepsilon/2})
\end{align*}
where $I_1$ and $I_2$ stand for the terms introduced at the proof of
Lemma~\ref{l6}, and the same arguments as therein lead to
$$
|\K_{sm,\varepsilon}f(x)|\leq C(n,p,\alpha) m^{n/2}\left( M(|f|^q)(x)
\right)^{1/q}
$$
for any $q\in(1,\infty).$
It remains to take the $L^{p,\omega}$-norms of the both sides above for $1<q<p$ and to
apply Lemma~\ref{l1} in order to get \eqref{Kpl1}.

The commutator estimate  follows analogously.
\end{pf}

Returning to the series expansions \eqref{EK}, we are in a
position now to complete the proof of Theorem~\ref{th1}.
First of all, note that
\begin{align*}
\sum_{m=1}^{\infty} \sum_{s=1}^{g_m}
\|b_{sm}\K_{sm,\varepsilon}f\|_{p,\omega}
\leq&\ C(n,p,\alpha,k) \|f\|_{p,\omega} \sum_{m=1}^{\infty} m^{-2l+n-2+n/2},\\
\sum_{m=1}^{\infty} \sum_{s=1}^{g_m}
	\|b_{sm}\C_{sm,\varepsilon}[a,k]f\|_{p,\omega}
\leq&\ C(n,p,\alpha,k) \|a\|_\ast \|f\|_{p,\omega} \sum_{m=1}^{\infty}
m^{-2l+n-2+n/2} \end{align*}
as it follows from \eqref{Y4}, \eqref{*} and Lemma~\ref{c1}.
Choosing $l>(3n-2)/4$ the series in \eqref{EK} result totally
convergent in $L^{p,\omega}(\R^n),$ uniformly in $\varepsilon>0,$ whence
$$
\|\K_{\varepsilon}f\|_{p,\omega} \leq C \|f\|_{p,\omega},\quad
\|\C_{\varepsilon}[a,k]f\|_{p,\omega} \leq C \|a\|_\ast\|f\|_{p,\omega}.
$$
Setting
$$
\K f(x) := \sum_{s,m} b_{sm}(x)\K_{sm}f(x),\qquad
\C[a,k]f(x) := \sum_{s,m} b_{sm}(x)\C_{sm}[a,k]f(x),
$$
we obtain as above
$$
\|\K f\|_{p,\omega} \leq C \|f\|_{p,\omega},\qquad
\|\C[a,k]f\|_{p,\omega} \leq C \|a\|_\ast\|f\|_{p,\omega}
$$
through \eqref{Y4}, \eqref{*} and Lemma~\ref{l6a}.

Finally, the total convergence in $L^{p,\omega}(\R^n)$ of the series
expansions
\eqref{EK}, uniformly in $\varepsilon>0,$ gives
\begin{align*}
&\lim_{\varepsilon\to0} \K_\varepsilon f(x) = \sum_{s,m} b_{sm}(x)
      \lim_{\varepsilon\to0} \K_{sm,\varepsilon}f(x)=
  \sum_{s,m} b_{sm}(x)
\K_{sm}f(x)=\K f(x),\\
&\lim_{\varepsilon\to0} \C_\varepsilon [a,k]f(x) = \sum_{s,m} b_{sm }(x)
      \lim_{\varepsilon\to0} \C_{sm,\varepsilon}[a,k]f(x)=\C[a,k]f(x)
\end{align*}
 and this completes the proof of Theorem~\ref{th1}.
\end{pf}

It is worth noting that singular integrals like \eqref{Kf} and \eqref{Cakf} appear in the representation formulas for the solutions of linear elliptic and parabolic
 partial differential equations. To make the obtained here results applicable to the study of regularizing properties of these operators we need of some additional local results.
\begin{crlr}\label{crl1}
Let $\Omega$ be a bounded domain in $\R^n$ and
$k(x;\xi)\colon\ \Omega\times (\R^n\setminus\{0\})\to\R$ be a variable kernel of
mixed homogeneity, $a\in  BMO(\Omega),$ $p\in (1,\infty)$ and
$\omega$ satisfies \eqref{Na1} and \eqref{Na2}. Then, for any $f\in
L^{p,\omega}(\Omega)$ and almost
all $x\in\Omega,$ the singular integrals
\begin{align*}
\K f(x)= & P.V. \int_{\Omega} k(x;x-y)f(y)dy\\
\C[a,k]f(x)= & P.V. \int_{\Omega} k(x;x-y)[a(y)-a(x)]f(y)dy
\end{align*}
are well defined in $L^{p,\omega}(\Omega)$ and
$$
\|\K f\|_{p,\omega;\Omega}   \leq C\|f\|_{p,\omega;\Omega},\quad
\|\C[a,k]f \|_{p,\omega;\Omega} \leq C\|a\|_\ast \|f\|_{p,\omega;\Omega}
$$
with $C=C(n,p,\alpha,\Omega,k).$
\end{crlr}
To obtain the above assertion it  is sufficient  to
extend $k(x;\cdot)$ and $f(\cdot)$ as zero outside $\Omega.$ One more necessary extension preserving the norm
 is that of $a$ in $BMO(\R^n)$
and we have it according to the results of Jones~\cite{PJ} and
Acquistapace~\cite{A} (see \cite{CFL} for details).

Another consequence of Theorem~\ref{th1} is the ``good behavior'' of the
commutator for $VMO$ functions $a.$
\begin{crlr}\label{crl2}
Suppose
$a\in  VMO$ with $VMO$-modulus $\gamma_a.$ Then, for each
$\varepsilon>0$
there exists  $r_0=r_0(\varepsilon,\gamma_a)>0$ such that for
any $\varrho\in(0,r_0)$ and any ellipsoid $\E_\varrho$ of radius $\varrho$ one has
\begin{equation}\label{vmo-omega}
\|\C[a,k]f \|_{p,\omega;\E_\varrho} \leq C \varepsilon \|f\|_{p,\omega;\E_\varrho}
\end{equation}
for all $f\in L^{p,\omega}(\E_\varrho).$
\end{crlr}
\begin{pf}
From the properties of the $VMO$  functions \cite[Theorem~1]{S} it follows that
for any $\varepsilon>0$  there exists $r_0=r_0(\varepsilon,\gamma_a)$  and
continuous and uniformly bounded function $g$  with modulus of continuity
$\omega_g(r_0)< \varepsilon/2$ such that $\|a-g\|_\ast<\varepsilon/2.$
Let $\E_\varrho$  be an ellipsoid  centered at $x_0$	and of radius $\varrho<r_0.$
Following \cite{CFL}  we   construct a function
$$
h(x)=\begin{cases}
g(x) & x\in\E_\varrho\\
g\left({x_0}_1+\varrho^{\alpha_1}\frac{x_1-{x_0}_1}{\rho(x-x_0)^{\alpha_1}},\ldots,
{x_0}_n+\varrho^{\alpha_n}\frac{x_n-{x_0}_n}{\rho(x-x_0)^{\alpha_n}}\right)
&x\in\E_\varrho^c
\end{cases}
$$
which is uniformly continuous in $\R^n.$  Whence the oscillation of $h$
in $\R^n$  is no greater than the
oscillation of $g$  in $\E_{r_0}.$   Then
\begin{align*}
\|\C[a,k]f \|_{p,\omega;\E_\varrho} & \leq \|\C[a-g,k]f \|_{p,\omega;\E_\varrho}+
\|\C[g,k]f\|_{p,\omega;\E_\varrho}\\
& \leq C \|a-g\|_\ast \|f\|_{p,\omega;\E_\varrho} + C\|h\|_\ast\|f\|_{p,\omega;\E_\varrho}\\
& \leq C\left(\|a-g\|_\ast +\omega_g(r_0) \right)\|f\|_{p,\omega;\E_\varrho} <
C\varepsilon \|f\|_{p,\omega;\E_\varrho}.
\end{align*}
\end{pf}

\end{document}